\documentclass{article}

\newtheorem{theorem}{Theorem}

\newtheorem{example}[theorem]{Example}

\input{tcilatex}
\begin{document}

\title{An Inversion Algorithm for the Cyclic Nonadiagonal Matrix}
\author{M. Ya\c{s}ar\thanks{%
e-mail: myasar@nigde.edu.tr} \& D. Bozkurt\thanks{%
e-mail: dbozkurt@selcuk.edu.tr} \and Department of Mathematics of Art and
Science Faculty, \and Ni\u{g}de University \& Department of Mathematics,
Science Faculty, \and Sel\c{c}uk University}
\maketitle

\begin{abstract}
In this paper, we compose a computational algorithm for the\linebreak
determinant and the inverse of the $n\times n$ cyclic nonadiagonal matrix.
The algorithm is suited for implementation using computer algebra systems
(CAS) such as Mathematica and Maple.
\end{abstract}

\bigskip AMS Clasification Number:

\section{Inroduction}

The $n\times n$ cyclic nonadiagonal matrix is as in the following form:

\begin{equation}
K=\left[ 
\begin{array}{cccccccccccc}
d_{1} & a_{1} & A_{1} & M_{1} & z_{1} & 0 & 0 & 0 & \cdots & 0 & B_{1} & 
b_{1} \\ 
b_{2} & d_{2} & a_{2} & A_{2} & M_{2} & z_{2} & 0 & 0 & \ddots & 0 & 0 & 
B_{2} \\ 
B_{3} & b_{3} & d_{3} & a_{3} & A_{3} & M_{3} & z_{3} & 0 & \ddots & \ddots
& \ddots & 0 \\ 
N_{4} & B_{4} & b_{4} & d_{4} & a_{4} & A_{4} & M_{4} & z_{4} & \ddots & 
\ddots & \ddots &  \\ 
R_{5} & N_{5} & B_{5} &  &  &  &  &  & \ddots & \ddots & \ddots & \vdots \\ 
0 & R_{6} & \ddots & \ddots & \ddots & \ddots & \ddots & \ddots & \ddots & 
& \ddots &  \\ 
&  & \ddots &  &  &  &  &  &  &  & \ddots & 0 \\ 
\vdots &  &  &  &  &  & \ddots & \ddots &  & \ddots &  & z_{n-4} \\ 
&  &  &  & \ddots &  &  &  &  &  &  & M_{n-3} \\ 
0 &  &  &  &  & \ddots &  &  &  & \ddots &  & A_{n-2} \\ 
A_{n-1} & 0 &  & \cdots &  & 0 & R_{n-1} & N_{n-1} & B_{n-1} & b_{n-1} & 
d_{n-1} & a_{n-1} \\ 
a_{n} & A_{n} & 0 & \cdots &  & \cdots & 0 & R_{n} & N_{n} & B_{n} & b_{n} & 
d_{n}%
\end{array}%
\right]  \label{1}
\end{equation}%
where $n\geq 8.$

This type of matrix appears in many areas such as engineering\linebreak
applications. The determinants and the inversions of these matrices
are\linebreak usually required. Many algorithms are composed by using the $%
LU $\linebreak factorization for the periodic tridiagonal, pentadiagonal and
the cyclic\linebreak pentadiagonal and heptadiagonal matrices\cite{1}-\cite%
{6}

A new recursive symbolic algorithm for inverting general periodic\linebreak
tridiagonal and anti-tridiagonal matrices are studied in \cite{2}. The
authors\linebreak compose a new symbolic algorithm for the inverses of the
periodic\linebreak pentadiagonal matrix and the periodic anti-pentadiagonal
matrix is obtained by using it in \cite{3}. With some restrictive
conditions, algorithms for the\linebreak inverses of the tridiagonal and
pentadiagonal matrices are given in \cite{5}. In \cite{4} it is\linebreak
presented that a new computational algorithm to evaluate the determinant of
the tridiagonal matrix with its cost. In \cite{6} an expression of
the\linebreak characteristic polynomial and eigenvectors for pentadiagonal
matrix is obtained and an algorithm to compute the determinant of the
pentadiagonal matrix is presented.

In this paper, we extend the work presented in \cite{1}. In the second
section we obtain the Doolittle $LU$ factorization of the cyclic
nonadiagonal matrix. Then by using the elements of the last six columns, the
elements of remaining $(n-6)$ columns are found and the inverse matrix of
the cyclic nonadiagonal matrix is obtained. In the last section a numerical
example is given.

\section{Main Result}

In this section, $t$ which is only a symbolic name is chosen as a parameter.
Then the determinant and the inverse of the cyclic nonadiagonal matrix $K$
in (\ref{1}) are computed. The $LU$ factorization of the matrix $K$ is as in
the following form:%
\[
L=\left[ 
\begin{array}{cccccccccc}
1 & 0 & 0 & 0 & 0 & 0 &  & \cdots & \cdots & 0 \\ 
f_{2} & 1 & 0 & 0 & \cdots & \cdots &  & \cdots & \cdots & 0 \\ 
g_{3} & f_{3} & 1 & 0 &  &  &  &  &  & \vdots \\ 
\alpha _{4} & g_{4} & f_{4} & 1 & \ddots &  &  &  &  & \vdots \\ 
\gamma _{5} & \alpha _{5} & g_{5} & f_{5} & 1 & \ddots &  &  &  & \vdots \\ 
0 & \ddots & \ddots & \ddots & \ddots & \ddots &  & \ddots &  & \vdots \\ 
\vdots &  & \ddots &  &  &  &  &  &  &  \\ 
0 & \ddots & 0 & \gamma _{n-2} & \alpha _{n-2} & g_{n-2} & f_{n-2} & 1 & 
\ddots & \vdots \\ 
k_{1} & k_{2} & \cdots & \cdots & \cdots &  & k_{n-3} & k_{n-2} & 1 & 0 \\ 
h_{1} & h_{2} & \cdots & \cdots & \cdots &  & h_{n-3} & h_{n-2} & h_{n-1} & 1%
\end{array}%
\right] 
\]%
and%
\[
U=\left[ 
\begin{array}{ccccccccccc}
c_{1} & e_{1} & P_{1} & T_{1} & z_{1} & 0 & \cdots & 0 & 0 & w_{1} & v_{1}
\\ 
0 & c_{2} & e_{2} & P_{2} & T_{2} & z_{2} & \cdots & 0 & 0 & w_{2} & v_{2}
\\ 
0 & \ddots & \ddots & \ddots & \ddots & \ddots & \ddots & \vdots &  & \vdots
& \vdots \\ 
&  & \ddots &  & \ddots &  & \ddots & z_{n-7} & 0 & w_{n-7} & v_{n-7} \\ 
\vdots & \ddots & \ddots &  & c_{n-6} & e_{n-6} & P_{n-6} & T_{n-6} & z_{n-6}
& w_{n-6} & v_{n-6} \\ 
\vdots &  & \ddots & \ddots &  & c_{n-5} & e_{n-5} & P_{n-5} & T_{n-5} & 
w_{n-5} & v_{n-5} \\ 
\vdots &  &  & \ddots & \ddots &  & c_{n-4} & e_{n-4} & P_{n-4} & w_{n-4} & 
v_{n-4} \\ 
\vdots &  &  &  & \ddots & \ddots &  & c_{n-3} & e_{n-3} & w_{n-3} & v_{n-3}
\\ 
\vdots &  &  &  &  & \ddots &  &  & c_{n-2} & w_{n-2} & v_{n-2} \\ 
&  &  &  &  &  & \ddots &  &  & c_{n-1} & v_{n-1} \\ 
0 & \cdots & \cdots & \cdots & \cdots & \cdots &  & \cdots &  & 0 & c_{n}%
\end{array}%
\right] 
\]%
where%
\begin{eqnarray*}
c_{i} &=&\left\{ 
\begin{array}{ll}
d_{1} & ,i=1 \\ 
d_{2}-f_{2}e_{1} & ,i=2 \\ 
d_{3}-g_{3}P_{1}-f_{3}e_{2} & ,i=3 \\ 
d_{4}-\alpha _{4}T_{1}-g_{4}P_{2}-f_{4}e_{3} & ,i=4 \\ 
d_{i}-\gamma _{i}z_{i-4}-\alpha _{i}T_{i-3}-g_{i}P_{i-2}-f_{i}e_{i-1} & 
,i=5,6,\ldots ,n-2 \\ 
d_{n-1}-\dsum\limits_{i=1}^{n-2}k_{i}w_{i} & ,i=n-1 \\ 
d_{n}-\dsum\limits_{i=1}^{n-1}h_{i}v_{i} & ,i=n%
\end{array}%
\right. \\
f_{i} &=&\left\{ 
\begin{array}{ll}
\dfrac{b_{2}}{c_{1}} & ,i=2 \\ 
\dfrac{1}{c_{2}}(b_{3}-g_{3}e_{1}) & ,i=3 \\ 
\dfrac{1}{c_{3}}(b_{4}-\alpha _{4}P_{1}-g_{4}e_{2}) & ,i=4 \\ 
\dfrac{1}{c_{i-1}}(b_{i}-\gamma _{i}T_{i-4}-\alpha _{i}P_{i-3}-g_{i}e_{i-2})
& ,i=5,6,\ldots ,n-2%
\end{array}%
\right. \\
g_{i} &=&\left\{ 
\begin{array}{ll}
\dfrac{B_{3}}{c_{1}} & ,i=3 \\ 
\dfrac{1}{c_{2}}(B_{4}-\alpha _{4}e_{1}) & ,i=4 \\ 
\dfrac{1}{c_{i-2}}(B_{i}-\gamma _{i}P_{i-4}-\alpha _{i}e_{i-3}) & 
,i=5,6,\ldots ,n-2%
\end{array}%
\right. \\
\alpha _{i} &=&\left\{ 
\begin{array}{ll}
\dfrac{N_{4}}{c_{1}} & ,i=4 \\ 
\dfrac{1}{c_{i-3}}(N_{i}-\gamma _{i}e_{i-4}) & ,i=5,6,\ldots ,n-2%
\end{array}%
\right. \\
e_{i} &=&\left\{ 
\begin{array}{ll}
a_{1} & ,i=1 \\ 
a_{2}-f_{2}P_{1} & ,i=2 \\ 
a_{3}-g_{3}T_{1}-f_{3}P_{2} & ,i=3 \\ 
a_{i}-\alpha _{i}z_{i-3}-g_{i}T_{i-2}-f_{i}P_{i-1} & ,i=4,5,\ldots ,n-3%
\end{array}%
\right. \\
P_{i} &=&\left\{ 
\begin{array}{ll}
A_{1} & ,i=1 \\ 
A_{2}-f_{2}T_{1} & ,i=2 \\ 
A_{i}-g_{i}z_{i-2}-f_{i}T_{i-1} & ,i=3,4,\ldots ,n-4%
\end{array}%
\right. \\
T_{i} &=&\left\{ 
\begin{array}{ll}
M_{1} & ,i=1 \\ 
M_{i}-f_{i}z_{i-1} & ,i=2,3,\ldots ,n-5%
\end{array}%
\right. \\
\gamma _{i} &=&\left\{ 
\begin{array}{cc}
\dfrac{R_{i}}{c_{i-4}} & ,i=5,6,\ldots ,n-2%
\end{array}%
\right.
\end{eqnarray*}%
\begin{eqnarray*}
k_{i} &=&\left\{ 
\begin{array}{ll}
\frac{A_{n-1}}{c_{1}} & ,i=1 \\ 
-\frac{1}{c_{2}}(k_{1}e_{1}) & ,i=2 \\ 
-\frac{1}{c_{3}}(k_{1}P_{1}+k_{2}e_{2}) & ,i=3 \\ 
-\frac{1}{c_{4}}(k_{1}T_{1}+k_{2}P_{2}+k_{3}e_{3}) & ,i=4 \\ 
-\frac{1}{c_{i}}(k_{i-4}z_{i-4}+k_{i-3}T_{i-3}+k_{i-2}P_{i-2}+k_{i-1}e_{i-1})
& ,i=5,\ldots ,n-6 \\ 
\frac{1}{c_{n-5}}%
(R_{n-1}-k_{n-9}z_{n-9}-k_{n-8}T_{n-8}-k_{n-7}P_{n-7}-k_{n-6}e_{n-6}) & 
,i=n-5 \\ 
\frac{1}{c_{n-4}}%
(N_{n-1}-k_{n-8}z_{n-8}-k_{n-7}T_{n-7}-k_{n-6}P_{n-6}-k_{n-5}e_{n-5}) & 
,i=n-4 \\ 
\frac{1}{c_{n-3}}%
(B_{n-1}-k_{n-7}z_{n-7}-k_{n-6}T_{n-6}-k_{n-5}P_{n-5}-k_{n-4}e_{n-4}) & 
,i=n-3 \\ 
\frac{1}{c_{n-2}}%
(b_{n-1}-k_{n-6}z_{n-6}-k_{n-5}T_{n-5}-k_{n-4}P_{n-4}-k_{n-3}e_{n-3}) & 
,i=n-2%
\end{array}%
\right. \\
w_{i} &=&\left\{ 
\begin{array}{ll}
B_{1} & ,i=1 \\ 
-f_{2}w_{1} & ,i=2 \\ 
-g_{3}w_{1}-f_{3}w_{2} & ,i=3 \\ 
-\alpha _{4}w_{1}-g_{4}w_{2}-f_{4}w_{3} & ,i=4 \\ 
-\gamma _{i}w_{i-4}-\alpha _{i}w_{i-3}-g_{i}w_{i-2}-f_{i}w_{i-1} & 
,i=5,\ldots ,n-6 \\ 
K_{n-5}-\gamma _{n-5}w_{n-9}-\alpha
_{n-5}w_{n-8}-g_{n-5}w_{n-7}-f_{n-5}w_{n-6} & ,i=n-5 \\ 
M_{n-4}-\gamma _{n-4}w_{n-8}-\alpha
_{n-4}w_{n-7}-g_{n-4}w_{n-6}-f_{n-4}w_{n-5} & ,i=n-4 \\ 
A_{n-3}-\gamma _{n-3}w_{n-7}-\alpha
_{n-3}w_{n-6}-g_{n-3}w_{n-5}-f_{n-3}w_{n-4} & ,i=n-3 \\ 
a_{n-2}-\gamma _{n-2}w_{n-6}-\alpha
_{n-2}w_{n-5}-g_{n-2}w_{n-4}-f_{n-2}w_{n-3} & ,i=n-2%
\end{array}%
\right. \\
h_{i} &=&\left\{ 
\begin{array}{ll}
\frac{a_{n}}{c_{1}} & ,i=1 \\ 
\frac{1}{c_{2}}(A_{n}-h_{1}e_{1}) & ,i=2 \\ 
-\frac{1}{c_{3}}(h_{1}P_{1}+h_{2}e_{2}) & ,i=3 \\ 
-\frac{1}{c_{4}}(h_{1}T_{1}+h_{2}P_{2}+h_{3}e_{3}) & ,i=4 \\ 
-\frac{1}{c_{i}}(h_{i-4}z_{i-4}+h_{i-3}T_{i-3}+h_{i-2}P_{i-2}+h_{i-1}e_{i-1})
& ,i=5,\ldots ,n-5 \\ 
\frac{1}{c_{n-4}}%
(R_{n}-h_{n-8}z_{n-8}-h_{n-7}T_{n-7}-h_{n-6}P_{n-6}-h_{n-5}e_{n-5}) & ,i=n-4
\\ 
\frac{1}{c_{n-3}}%
(N_{n}-h_{n-7}z_{n-7}-h_{n-6}T_{n-6}-h_{n-5}P_{n-5}-h_{n-4}e_{n-4}) & ,i=n-3
\\ 
\frac{1}{c_{n-2}}%
(B_{n}-h_{n-6}z_{n-6}-h_{n-5}T_{n-5}-h_{n-4}P_{n-4}-h_{n-3}e_{n-3}) & ,i=n-2
\\ 
\frac{1}{c_{n-1}}(b_{n}-\dsum\limits_{i=1}^{n-2}h_{i}w_{i}) & ,i=n-1%
\end{array}%
\right. \\
v_{i} &=&\left\{ 
\begin{array}{ll}
b_{1} & ,i=1 \\ 
B_{2}-f_{2}v_{1} & ,i=2 \\ 
-g_{3}v_{1}-f_{3}v_{2} & ,i=3 \\ 
-\alpha _{4}v_{1}-g_{4}v_{2}-f_{4}v_{3} & ,i=4 \\ 
-\gamma _{i}v_{i-4}-\alpha _{i}v_{i-3}-g_{i}v_{i-2}-f_{i}v_{i-1} & 
,i=5,\ldots ,n-5 \\ 
K_{n-4}-\gamma _{n-4}v_{n-8}-\alpha
_{n-4}v_{n-7}-g_{n-4}v_{n-6}-f_{n-4}v_{n-5} & ,i=n-4 \\ 
M_{n-3}-\gamma _{n-3}v_{n-7}-\alpha
_{n-3}v_{n-6}-g_{n-3}v_{n-5}-f_{n-3}v_{n-4} & ,i=n-3 \\ 
A_{n-2}-\gamma _{n-2}v_{n-6}-\alpha
_{n-2}v_{n-5}-g_{n-2}v_{n-4}-f_{n-2}v_{n-3} & ,i=n-2 \\ 
a_{n-1}-\dsum\limits_{i=1}^{n-2}k_{i}v_{i} & ,i=n-1%
\end{array}%
\right.
\end{eqnarray*}

\newpage

The determinant of the cyclic nonadiagonal matrix $K$ is computed as in the
following: 
\[
\det (K)=\dprod\limits_{i=1}^{n}c_{i} 
\]

When $K$ is nonsingular, let%
\[
K^{-1}=(S_{ij})_{1\leq i,j\leq n}=(C_{1},C_{2},\ldots ,C_{r},\ldots ,C_{n}) 
\]%
and $C_{i}$ is the $i$th column of $K^{-1}$ where $1\leq i\leq n$. Notice
that%
\[
C_{r}=(S_{1,r},S_{2,r},...,S_{n,r})^{T} 
\]%
for $r=1,2,...,n.$ We can write $C_{r}$ as follows:%
\begin{equation}
C_{r}=(C_{1},C_{2},\ldots ,C_{r},\ldots ,C_{n})E_{r}  \label{a}
\end{equation}%
where $E_{r}=(\delta _{1r},\delta _{2r},\ldots ,\delta _{rr},\ldots ,\delta
_{1n})^{T},r=1,2,\ldots ,n$ ($\delta _{ij}$ is the Kronecker\linebreak
symbol). By using (\ref{a}) we obtain 
\begin{equation}
KC_{i}=E_{i}  \label{b}
\end{equation}%
for $i=n,n-1,n-2,n-3,n-4,n-5.$

Now, an algorithm can be composed by using the last six columns of $K^{-1}.$
By using the $LU$ factorization and (\ref{b}), the components of the last
six columns are computed as in the following:%
\begin{equation}
\begin{array}{l}
S_{n,n}=\frac{1}{c_{n}} \\ 
S_{n-1,n}=-\frac{1}{c_{n-1}}(v_{n-1}S_{n,n}) \\ 
S_{n,n-1}=-\frac{h_{n-1}}{c_{n}} \\ 
S_{n-1,n-1}=\frac{1}{c_{n-1}}(1-v_{n-1}S_{n,n-1}) \\ 
S_{n,n-2}=\frac{1}{c_{n}}(-h_{n-2}+h_{n-1}k_{n-2}) \\ 
S_{n-1,n-2}=\frac{1}{c_{n-1}}(-k_{n-2}-v_{n-1}S_{n,n-2}) \\ 
S_{n-2,n-2}=\frac{1}{c_{n-2}}(1-w_{n-2}S_{n-1,n-2}-v_{n-2}S_{n,n-2}) \\ 
S_{n,n-3}=\frac{1}{c_{n}}%
(-h_{n-3}+h_{n-2}f_{n-2}+h_{n-1}k_{n-3}-h_{n-1}k_{n-2}f_{n-2}) \\ 
S_{n-1,n-3}=\frac{1}{c_{n-1}}(-k_{n-3}+k_{n-2}f_{n-2}-v_{n-1}S_{n,n-3}) \\ 
S_{n-2,n-3}=\frac{1}{c_{n-2}}(-f_{n-2}-w_{n-2}S_{n-1,n-3}-v_{n-2}S_{n,n-3})
\\ 
S_{n-3,n-3}=\frac{1}{c_{n-3}}%
(1-e_{n-3}S_{n-2,n-3}-w_{n-3}S_{n-1,n-3}-v_{n-3}S_{n,n-3}) \\ 
\begin{array}{l}
S_{n,n-4}=\frac{1}{c_{n}}%
(-h_{n-4}+h_{n-3}f_{n-3}+h_{n-2}g_{n-2}-h_{n-2}f_{n-2}f_{n-3} \\ 
+h_{n-1}k_{n-4}-h_{n-1}k_{n-3}f_{n-3}-h_{n-1}k_{n-2}g_{n-2}+h_{n-1}k_{n-2}f_{n-2}f_{n-3})%
\end{array}
\\ 
S_{n-1,n-4}=\frac{1}{c_{n-1}}%
(-k_{n-4}+k_{n-3}f_{n-3}+k_{n-2}g_{n-2}-k_{n-2}f_{n-2}f_{n-3} \\ 
-v_{n-1}S_{n,n-4})%
\end{array}
\label{3}
\end{equation}%
\begin{equation}
\begin{array}{l}
S_{n-2,n-4}=\frac{1}{c_{n-2}}%
(-g_{n-2}+f_{n-2}f_{n-3}-w_{n-2}S_{n-1,n-4}-v_{n-2}S_{n,n-4}) \\ 
\begin{array}{l}
S_{n-3,n-4}=\frac{1}{c_{n-3}}(-f_{n-3}-e_{n-3}S_{n-2,n-4} \\ 
-w_{n-3}S_{n-1,n-4}-v_{n-3}S_{n,n-4})%
\end{array}
\\ 
\begin{array}{l}
S_{n-4,n-4}=\frac{1}{c_{n-4}}(1-e_{n-4}S_{n-3,n-4}-P_{n-4}S_{n-2,n-4} \\ 
-w_{n-4}S_{n-1,n-4}-v_{n-4}S_{n,n-4})%
\end{array}
\\ 
\begin{array}{l}
S_{n,n-5}=\frac{1}{c_{n}}%
(-h_{n-5}+h_{n-4}f_{n-4}+h_{n-3}g_{n-3}-h_{n-3}f_{n-3}f_{n-4} \\ 
+h_{n-2}\alpha
_{n-2}-h_{n-2}g_{n-2}f_{n-4}-h_{n-2}f_{n-2}g_{n-3}+h_{n-2}f_{n-2}f_{n-3}f_{-4}
\\ 
+h_{n-1}k_{n-5}-h_{n-1}k_{n-4}f_{n-4}-h_{n-1}k_{n-3}g_{n-3}+h_{n-1}k_{n-3}f_{n-3}f_{n-4}
\\ 
-h_{n-1}k_{n-2}\alpha
_{n-2}+h_{n-1}k_{n-2}g_{n-2}f_{n-4}+h_{n-1}k_{n-2}f_{n-2}g_{n-3} \\ 
-h_{n-1}k_{n-2}f_{n-2}f_{n-3}f_{n-4})%
\end{array}
\\ 
\begin{array}{l}
S_{n-1,n-5}=\frac{1}{c_{n-1}}%
(-k_{n-5}+k_{n-4}f_{n-4}+k_{n-3}g_{n-3}-k_{n-3}f_{n-3}f_{n-4} \\ 
+k_{n-2}\alpha
_{n-2}-k_{n-2}g_{n-2}f_{n-4}-k_{n-2}f_{n-2}g_{n-3}+k_{n-2}f_{n-2}f_{n-3}f_{n-4}
\\ 
-v_{n-1}S_{n,n-5})%
\end{array}
\\ 
\begin{array}{l}
S_{n-2,n-5}=\frac{1}{c_{n-2}}(-\alpha
_{n-2}+g_{n-2}f_{n-4}+f_{n-2}g_{n-3}-f_{n-2}f_{n-3}f_{n-4} \\ 
-w_{n-2}S_{n-1,n-5}-v_{n-2}S_{n,n-5})%
\end{array}
\\ 
\begin{array}{l}
S_{n-3,n-5}=\frac{1}{c_{n-3}}(-g_{n-3}+f_{n-3}f_{n-4}-e_{n-3}S_{n-2,n-5} \\ 
-w_{n-3}S_{n-1,n-5}-v_{n-3}S_{n,n-5})%
\end{array}
\\ 
\begin{array}{l}
S_{n-4,n-5}=\frac{1}{c_{n-4}}(-f_{n-4}-e_{n-4}S_{n-3,n-5}-P_{n-4}S_{n-2,n-5}
\\ 
-w_{n-4}S_{n-1,n-5}-v_{n-4}S_{n,n-5})%
\end{array}
\\ 
\begin{array}{l}
S_{n-5,n-5}=\frac{1}{c_{n-5}}(1-e_{n-5}S_{n-4,n-5}-P_{n-5}S_{n-3,n-5} \\ 
-T_{n-5}S_{n-2,n-5}-w_{n-5}S_{n-1,n-5}-v_{n-5}S_{n,n-5})%
\end{array}%
\end{array}
\label{4}
\end{equation}%
and for $j=n,n-1$%
\begin{equation}
S_{n-2,j}=-\frac{1}{c_{n-2}}(w_{n-2}S_{n-1,j}+v_{n-2}S_{n,j})  \label{5}
\end{equation}%
for $j=n,n-1,n-2$%
\begin{equation}
S_{n-3,j}=-\frac{1}{c_{n-3}}%
(e_{n-3}S_{n-2,j}+w_{n-3}S_{n-1,j}+v_{n-3}S_{n,j})  \label{6}
\end{equation}%
for $j=n,n-1,n-2,n-3$%
\begin{equation}
S_{n-4,j}=-\frac{1}{c_{n-4}}%
(e_{n-4}S_{n-3,j}+P_{n-4}S_{n-2,j}+w_{n-4}S_{n-1,j}+v_{n-4}S_{n,j})
\label{7}
\end{equation}%
for $j=n,n-1,n-2,n-3,n-4$%
\begin{equation}
S_{n-5,j}=-\frac{1}{c_{n-5}}%
(e_{n-5}S_{n-4,j}+P_{n-5}S_{n-3,j}+T_{n-5}S_{n-2,j}+w_{n-5}S_{n-1,j}+v_{n-5}S_{n,j})
\label{8}
\end{equation}%
for $j=n,n-1,n-2,n-3,n-4,n-5$ and $i=n-6,n-7,...,1$%
\begin{equation}
S_{i,j}=-\frac{1}{c_{i}}%
(e_{i}S_{i+1,j}+P_{i}S_{i+2,j}+T_{i}S_{i+3,j}+z_{i}S_{i+4,j}+w_{i}S_{n-1,j}+v_{i}S_{n,j})
\label{9}
\end{equation}

The elements of the remaining $(n-6)$ columns can be calculated using the
fact that $K^{-1}K=I_{n}$ where $I_{n}$ is an $n\times n$ identity matrix.
Then

\bigskip 
\begin{eqnarray}
C_{n-6} &=&\frac{1}{z_{n-6}}%
(E_{n-2}-M_{n-5}C_{n-5}-A_{n-4}C_{n-4}-a_{n-3}C_{n-3}  \label{10} \\
&&-d_{n-2}C_{n-2}-b_{n-1}C_{n-1}-B_{n}C_{n})  \nonumber
\end{eqnarray}%
\begin{eqnarray}
C_{n-7} &=&\frac{1}{z_{n-7}}%
(E_{n-3}-M_{n-6}C_{n-6}-A_{n-5}C_{n-5}-a_{n-4}C_{n-4}  \label{11} \\
&&-d_{n-3}C_{n-3}-b_{n-2}C_{n-2}-B_{n-1}C_{n-1}-N_{n}C_{n})  \nonumber
\end{eqnarray}

\begin{eqnarray}
C_{j} &=&\frac{1}{z_{j}}%
(E_{j+4}-M_{j+1}C_{j+1}-A_{j+2}C_{j+2}-a_{j+3}C_{j+3}-d_{j+4}C_{j+4}
\label{12} \\
&&-b_{j+5}C_{j+5}-B_{j+6}C_{j+6}-N_{j+7}C_{j+7}-R_{j+8}C_{j+8})  \nonumber
\end{eqnarray}%
where $j=n-8,n-9,\ldots ,1$ and $z_{i}\neq 0$ for $i=1,2,\ldots ,n-6.$

Thus, the algorithm for the $n\times n$ cyclic nonadiagonal matrix is given
as

\textbf{Input: }$n$ is the order and $%
d_{i},a_{i},A_{i},M_{i},z_{i},b_{i},B_{i},N_{i},R_{i}$ are the entries of
the cyclic nonadiagonal matrix.

\textbf{Output: }Inverse matrix $K^{-1}=(S_{ij})_{1\leq i,j\leq n}.$

\textbf{Step1:}If $z_{i}=0$ for any $i=1,2,\ldots ,n-6$ set $z_{i}=t.$

\textbf{Step2:} If $R_{i}=0$ for any $i=6,7,\ldots ,n$ set $R_{i}=t.$

\textbf{Step3:} Set $c_{1}=d_{1},$ if $c_{1}=0$ then $c_{1}=t,~$%
\begin{eqnarray*}
f_{2} &=&\dfrac{b_{2}}{c_{1}} \\
g_{3} &=&\dfrac{B_{3}}{c_{1}} \\
e_{1} &=&a_{1} \\
P_{1} &=&A_{1} \\
\alpha _{4} &=&\frac{N_{4}}{c_{1}} \\
T_{1} &=&M_{1} \\
k_{1} &=&\frac{A_{n-1}}{c_{1}} \\
w_{1} &=&B_{1} \\
h_{1} &=&\frac{a_{n}}{c_{1}} \\
v_{1} &=&b_{1}
\end{eqnarray*}%
$c_{2}=d_{2}-f_{2}e_{1},$ if $c_{2}=0$ then $c_{2}=t,$%
\begin{eqnarray*}
f_{3} &=&\dfrac{1}{c_{2}}(b_{3}-g_{3}e_{1}) \\
g_{4} &=&\frac{1}{c_{2}}(B_{4}-\alpha _{4}e_{1}) \\
e_{2} &=&a_{2}-f_{2}P_{1} \\
P_{2} &=&A_{2}-f_{2}T_{1} \\
T_{2} &=&M_{2}-f_{2}z_{1} \\
k_{2} &=&-\frac{1}{c_{2}}(k_{1}e_{1}) \\
w_{2} &=&-f_{2}w_{1} \\
h_{2} &=&\frac{1}{c_{2}}(A_{n}-h_{1}e_{1}) \\
v_{2} &=&B_{2}-f_{2}v_{1}
\end{eqnarray*}%
$c_{3}=d_{3}-g_{3}P_{1}-f_{3}e_{2},$ if $c_{3}=0$ then $c_{3}=t$%
\begin{eqnarray*}
f_{4} &=&\frac{1}{c_{3}}(b_{4}-\alpha _{4}P_{1}-g_{4}e_{2}) \\
e_{3} &=&a_{3}-g_{3}T_{1}-f_{3}P_{2} \\
P_{3} &=&A_{3}-g_{3}z_{1}-f_{3}T_{2} \\
T_{3} &=&M_{3}-f_{3}z_{2} \\
k_{3} &=&-\frac{1}{c_{3}}(k_{1}P_{1}+k_{2}e_{2}) \\
w_{3} &=&-g_{3}w_{1}-f_{3}w_{2} \\
h_{3} &=&-\frac{1}{c_{3}}(h_{1}P_{1}+h_{2}e_{2}) \\
v_{3} &=&-g_{3}v_{1}-f_{3}v_{2}
\end{eqnarray*}%
$c_{4}=d_{4}-\alpha _{4}T_{1}-g_{4}P_{2}-f_{4}e_{3}$ if $c_{4}=0$ then $%
c_{4}=t$%
\begin{eqnarray*}
e_{4} &=&a_{4}-\alpha _{4}z_{1}-g_{4}T_{2}-f_{4}P_{3} \\
P_{4} &=&A_{4}-g_{4}z_{2}-f_{4}T_{3} \\
T_{4} &=&M_{4}-f_{4}z_{3} \\
k_{4} &=&-\frac{1}{c_{4}}(k_{1}T_{1}+k_{2}P_{2}+k_{3}e_{3}) \\
w_{4} &=&-\alpha _{4}w_{1}-g_{4}w_{2}-f_{4}w_{3} \\
h_{4} &=&-\frac{1}{c_{4}}(h_{1}T_{1}+h_{2}P_{2}+h_{3}e_{3}) \\
v_{4} &=&-\alpha _{4}v_{1}-g_{4}v_{2}-f_{4}v_{3}
\end{eqnarray*}

\bigskip \textbf{Step4:} For $i=5,6,\ldots ,n-3$ do

$\gamma _{i}=\dfrac{R_{i}}{c_{i-4}}$

$\alpha _{i}=\dfrac{1}{c_{i-3}}(N_{i}-\gamma _{i}e_{i-4})$

$g_{i}=\dfrac{1}{c_{i-2}}(B_{i}-\gamma _{i}P_{i-4}-\alpha _{i}e_{i-3})$

$f_{i}=\dfrac{1}{c_{i-1}}(b_{i}-\gamma _{i}T_{i-4}-\alpha
_{i}P_{i-3}-g_{i}e_{i-2})$

$e_{i}=a_{i}-\alpha _{i}z_{i-3}-g_{i}T_{i-2}-f_{i}P_{i-1}$

$c_{i}=d_{i}-\gamma _{i}z_{i-4}-\alpha _{i}T_{i-3}-g_{i}P_{i-2}-f_{i}e_{i-1}$
if $c_{i}=0$ then $c_{i}=t$

and set

$\gamma _{n-2}=\frac{R_{n-2}}{c_{n-6}}$

$\alpha _{n-2}=\frac{1}{c_{n-5}}(N_{n-2}-\gamma _{n-2}e_{n-6})$

$g_{n-2}=\dfrac{1}{c_{n-4}}(B_{n-2}-\gamma _{n-2}P_{n-6}-\alpha
_{n-2}e_{n-5})$

$f_{n-2}=\dfrac{1}{c_{n-3}}(b_{n-2}-\gamma _{n-2}T_{n-6}-\alpha
_{n-2}P_{n-5}-g_{n-2}e_{n-4})$

$c_{n-2}=d_{n-2}-\gamma _{n-2}z_{n-6}-\alpha
_{n-2}T_{n-5}-g_{n-2}P_{n-4}-f_{n-2}e_{n-3},$ if $c_{n-2}=0$ then $%
c_{n-2}=t. $

\textbf{Step5:} $i=5,6\ldots n-6$

$k_{i}=-\frac{1}{c_{i}}%
(k_{i-4}z_{i-4}+k_{i-3}T_{i-3}+k_{i-2}P_{i-2}+k_{i-1}e_{i-1})$

$w_{i}=-\gamma _{i}w_{i-4}-\alpha _{i}w_{i-3}-g_{i}w_{i-2}-f_{i}w_{i-1}$

$k_{n-5}=\frac{1}{c_{n-5}}%
(R_{n-1}-k_{n-9}z_{n-9}-k_{n-8}T_{n-8}-k_{n-7}P_{n-7}-k_{n-6}e_{n-6})$

$k_{n-4}=\frac{1}{c_{n-4}}%
(N_{n-1}-k_{n-8}z_{n-8}-k_{n-7}T_{n-7}-k_{n-6}P_{n-6}-k_{n-5}e_{n-5})$

$k_{n-3}=\frac{1}{c_{n-3}}%
(B_{n-1}-k_{n-7}z_{n-7}-k_{n-6}T_{n-6}-k_{n-5}P_{n-5}-k_{n-4}e_{n-4})$

$k_{n-2}=\frac{1}{c_{n-2}}%
(b_{n-1}-k_{n-6}z_{n-6}-k_{n-5}T_{n-5}-k_{n-4}P_{n-4}-k_{n-3}e_{n-3})$

$w_{n-5}=K_{n-5}-\gamma _{n-5}w_{n-9}-\alpha
_{n-5}w_{n-8}-g_{n-5}w_{n-7}-f_{n-5}w_{n-6}$

$w_{n-4}=M_{n-4}-\gamma _{n-4}w_{n-8}-\alpha
_{n-4}w_{n-7}-g_{n-4}w_{n-6}-f_{n-4}w_{n-5}$

$w_{n-3}=A_{n-3}-\gamma _{n-3}w_{n-7}-\alpha
_{n-3}w_{n-6}-g_{n-3}w_{n-5}-f_{n-3}w_{n-4}$

$w_{n-2}=a_{n-2}-\gamma _{n-2}w_{n-6}-\alpha
_{n-2}w_{n-5}-g_{n-2}w_{n-4}-f_{n-2}w_{n-3}$

$c_{n-1}=d_{n-1}-\dsum\limits_{i=1}^{n-2}k_{i}w_{i}$ if $c_{n-1}=0$ then $%
c_{n-1}=t.$

\textbf{Step6:} $i=5,6,\ldots ,n-5$

$T_{i}=M_{i}-f_{i}z_{i-1}$

$P_{i}=A_{i}-g_{i}z_{i-2}-f_{i}T_{i-1}$

$h_{i}=-\frac{1}{c_{i}}%
(h_{i-4}z_{i-4}+h_{i-3}T_{i-3}+h_{i-2}P_{i-2}+h_{i-1}e_{i-1})$

$v_{i}=-\gamma _{i}v_{i-4}-\alpha _{i}v_{i-3}-g_{i}v_{i-2}-f_{i}v_{i-1}$

$P_{n-4}=A_{n-4}-g_{n-4}z_{n-6}-f_{n-4}T_{n-5}$

$h_{n-4}=\frac{1}{c_{n-4}}%
(R_{n}-h_{n-8}z_{n-8}-h_{n-7}T_{n-7}-h_{n-6}P_{n-6}-h_{n-5}e_{n-5})$

$h_{n-3}=\frac{1}{c_{n-3}}%
(N_{n}-h_{n-7}z_{n-7}-h_{n-6}T_{n-6}-h_{n-5}P_{n-5}-h_{n-4}e_{n-4})$

$h_{n-2}=\frac{1}{c_{n-2}}%
(B_{n}-h_{n-6}z_{n-6}-h_{n-5}T_{n-5}-h_{n-4}P_{n-4}-h_{n-3}e_{n-3})$

$h_{n-1}=\frac{1}{c_{n-1}}(b_{n}-\dsum\limits_{i=1}^{n-2}h_{i}w_{i})$

$v_{n-4}=K_{n-4}-\gamma _{n-4}v_{n-8}-\alpha
_{n-4}v_{n-7}-g_{n-4}v_{n-6}-f_{n-4}v_{n-5}$

$v_{n-3}=M_{n-3}-\gamma _{n-3}v_{n-7}-\alpha
_{n-3}v_{n-6}-g_{n-3}v_{n-5}-f_{n-3}v_{n-4}$

$v_{n-2}=A_{n-2}-\gamma _{n-2}v_{n-6}-\alpha
_{n-2}v_{n-5}-g_{n-2}v_{n-4}-f_{n-2}v_{n-3}$

$v_{n-1}=a_{n-1}-\dsum\limits_{i=1}^{n-2}k_{i}v_{i}$

$c_{n}=d_{n}-\dsum\limits_{i=1}^{n-1}h_{i}v_{i}$ if $c_{n}=0$ then $c_{n}=t.$

\textbf{Step7:} Compute $\det (K)=\dprod\limits_{i=1}^{n}c_{i}.$ If $\det
(K)=0,$ then OUTPUT(Singular Matrix) stop.

\textbf{Step8:} For $i=1,2,\ldots ,n$ compute and simplify the components $%
S_{i,n},S_{i,n-1},$\linebreak $S_{i,n-2},S_{i,n-3},S_{i,n-4},S_{i,n-5}$ of
the columns $C_{j}$ where $j=n,n-1,n-2,n-3,n-4,n-5$ by using (\ref{3},\ref{4}%
,\ref{5},\ref{6},\ref{7},\ref{8},\ref{9}).

\textbf{Step9:} Compute the components of the columns $C_{n-6}$ and $C_{n-7}$
by using (\ref{10},\ref{11}), then for $j=n-8,n-9,\ldots ,1$ and $%
i=1,2,\ldots ,n$ compute and simplify the components $S_{ij}$ by using (\ref%
{12}).

\textbf{Step10:} For $i,j=1,2,\ldots ,n$ substitute the actual value $t=0$
in all $S_{i,j}.$

Let $R$ be an $n\times n$ matrix as in the following form:%
\[
R=\left[ 
\begin{array}{ccccc}
0 & \cdots & \cdots & 0 & 1 \\ 
\vdots &  & {\mathinner{\mkern2mu\raise1pt\hbox{.}\mkern2mu
\raise4pt\hbox{.}\mkern2mu\raise7pt\hbox{.}\mkern1mu}} & 1 & 0 \\ 
\vdots & {\mathinner{\mkern2mu\raise1pt\hbox{.}\mkern2mu
\raise4pt\hbox{.}\mkern2mu\raise7pt\hbox{.}\mkern1mu}} & {%
\mathinner{\mkern2mu\raise1pt\hbox{.}\mkern2mu
\raise4pt\hbox{.}\mkern2mu\raise7pt\hbox{.}\mkern1mu}} & {%
\mathinner{\mkern2mu\raise1pt\hbox{.}\mkern2mu
\raise4pt\hbox{.}\mkern2mu\raise7pt\hbox{.}\mkern1mu}} & \vdots \\ 
0 & {\mathinner{\mkern2mu\raise1pt\hbox{.}\mkern2mu
\raise4pt\hbox{.}\mkern2mu\raise7pt\hbox{.}\mkern1mu}} & {%
\mathinner{\mkern2mu\raise1pt\hbox{.}\mkern2mu
\raise4pt\hbox{.}\mkern2mu\raise7pt\hbox{.}\mkern1mu}} &  & \vdots \\ 
1 & 0 & \cdots & \cdots & 0%
\end{array}%
\right] 
\]%
It is clear that $R$ is a nonsingular and its inverse matrix is itself. Let $%
Y$ be a cyclic anti-nonadiagonal matrix. Since there is the following
relation between the cyclic nonadiagonal and the cyclic anti-nonadiagonal
matrices%
\[
Y=KR, 
\]%
the inverse matrix of $Y$ is obtained as%
\[
Y^{-1}=RK^{-1}. 
\]%
\newpage

\section{\protect\bigskip Numerical Example}

\begin{example}
Consider the matrix D as in the following%
\[
K=\left[ 
\begin{array}{rrrrrrrrrrrr}
1 & 1 & -1 & 2 & 1 & 0 & 0 & 0 & 0 & 0 & 1 & 1 \\ 
-2 & 2 & 1 & 1 & -1 & 2 & 0 & 0 & 0 & 0 & 0 & 1 \\ 
1 & 1 & 2 & 1 & 1 & -1 & 1 & 0 & 0 & 0 & 0 & 0 \\ 
1 & 2 & 1 & -1 & 1 & 1 & 1 & 1 & 0 & 0 & 0 & 0 \\ 
-1 & 1 & 1 & 2 & 1 & -1 & -1 & 1 & 1 & 0 & 0 & 0 \\ 
0 & 1 & -1 & 1 & 1 & 1 & 2 & -1 & 1 & 2 & 0 & 0 \\ 
0 & 0 & -1 & -1 & 2 & 1 & -1 & 2 & 1 & 1 & 1 & 0 \\ 
0 & 0 & 0 & 1 & 1 & -1 & 1 & 1 & 2 & 3 & -1 & 2 \\ 
0 & 0 & 0 & 0 & -2 & 1 & 1 & -1 & 2 & 1 & 1 & 1 \\ 
0 & 0 & 0 & 0 & 0 & 1 & -1 & 3 & 1 & -1 & -2 & 1 \\ 
1 & 0 & 0 & 0 & 0 & 0 & 1 & 2 & 3 & 1 & -1 & 1 \\ 
1 & 1 & 0 & 0 & 0 & 0 & 0 & 1 & 1 & 1 & -1 & 1%
\end{array}%
\right] . 
\]%
We apply the algorithm to it and we have
\end{example}

\begin{itemize}
\item $%
\begin{array}{ll}
(c_{1},c_{2},c_{3},c_{4})=(1,4,3,-\frac{7}{2}) & (f_{2},f_{3},f_{4})=(-2,0,%
\frac{3}{4}) \\ 
(g_{3},g_{4})=(1,\frac{1}{4}) & \alpha _{4}=1 \\ 
(e_{1},e_{2},e_{3},e_{4})=(1,-1,-1,-\frac{1}{4}) & 
(T_{1},T_{2},T_{3},T_{4})=(2,1,-1,\frac{1}{4}) \\ 
(k_{1},k_{2},k_{3},k_{4})=(1,-\frac{1}{4},\frac{1}{4},\frac{1}{7}) & 
(w_{1},w_{2},w_{3},w_{4})=(1,2,-1,-\frac{3}{4}) \\ 
(h_{1},h_{2},h_{3}h_{4})=(1,0,\frac{1}{3},\frac{10}{21}) & 
(P_{1},P_{2},P_{3},P_{4})=(-1,5,0,\frac{5}{4}) \\ 
(v_{1},v_{2},v_{3},v_{4})=(1,3,-1,-1) & 
\end{array}%
$

\item $%
\begin{array}{l}
(\gamma _{5},\gamma _{6},\gamma _{7},\gamma _{8},\gamma _{9},\gamma
_{10})=(-1,\frac{1}{4},-\frac{1}{3},-\frac{2}{7},-\frac{42}{29},\frac{58}{45}%
) \\ 
(\alpha _{5},\alpha _{6},\alpha _{7},\alpha _{8},\alpha _{9},\alpha _{10})=(%
\frac{1}{2},-\frac{1}{4},\frac{8}{21},\frac{39}{58},-\frac{46}{45},\frac{52}{%
75}) \\ 
(g_{5},g_{6},g_{7},g_{8},g_{9},g_{10})=(\frac{1}{6},\frac{1}{7},\frac{44}{29}%
,\frac{11}{45},-\frac{53}{150},\frac{116}{77}) \\ 
(f_{5},f_{6},f_{7},f_{8},f_{9},f_{10})=(-\frac{10}{21},\frac{33}{58},\frac{8%
}{3},-\frac{49}{300},\frac{38}{77},\frac{13}{218}) \\ 
(e_{5},e_{6},e_{7},e_{8},e_{9})=(-\frac{26}{21},\frac{163}{58},\frac{14}{3},%
\frac{19}{20},\frac{48}{77}) \\ 
(c_{5},c_{6},c_{7},c_{8},c_{9},c_{10})=(\frac{29}{21},\frac{45}{58},-\frac{20%
}{3},\frac{77}{50},\frac{218}{77},-\frac{1088}{327})%
\end{array}%
$

\item $%
\begin{array}{ll}
(k_{5},k_{6},k_{7},k_{8},k_{9},k_{10})=(-\frac{15}{29},-\frac{4}{45},-\frac{%
19}{300},\frac{137}{77},\frac{271}{436},\frac{1761}{2176}) & c_{11}=\frac{511%
}{544} \\ 
(w_{5},w_{6},w_{7},w_{8},w_{9},w_{10})=(-\frac{4}{21},-\frac{31}{58},\frac{8%
}{3},-\frac{13}{25},\frac{106}{77},-\frac{268}{109}) & 
\end{array}%
$

\item $%
\begin{array}{l}
(T_{5},T_{6},T_{7})=(\frac{31}{21},\frac{25}{58},-\frac{13}{3}) \\ 
(P_{5},P_{6},P_{7},P_{8})=(-\frac{22}{21},-\frac{115}{58},-\frac{5}{3},\frac{%
541}{300}) \\ 
(h_{5},h_{6},h_{7},h_{8},h_{9},h_{10},h_{11})=(-\frac{37}{58},-\frac{61}{45}%
,-\frac{121}{300},\frac{3}{7},\frac{44}{109},-\frac{307}{1088},-\frac{1146}{%
511}) \\ 
(v_{5},v_{6},v_{7},v_{8},v_{9},v_{10},v_{11})=(-\frac{17}{21},-\frac{23}{58},%
\frac{7}{3},\frac{821}{300},-\frac{85}{77},-\frac{2723}{654},\frac{101}{4352}%
) \\ 
c_{12}=-\frac{4715}{4088}%
\end{array}%
$

\item $\det (K)=4715$

\item $K^{-1}=\left[ 
\begin{array}{rrrrrr}
\frac{231}{4715} & \frac{199}{943} & \frac{3154}{4715} & -\frac{3181}{4715}
& -\frac{2187}{4715} & \frac{142}{4715} \\ 
\frac{79}{943} & -\frac{170}{943} & -\frac{591}{943} & \frac{643}{943} & 
\frac{440}{943} & -\frac{29}{943} \\ 
-\frac{1172}{4715} & \frac{15}{943} & \frac{3062}{4715} & -\frac{1088}{4715}
& -\frac{416}{4715} & \frac{96}{4715} \\ 
\frac{33}{205} & \frac{46}{41} & \frac{187}{205} & -\frac{308}{205} & -\frac{%
166}{205} & -\frac{9}{205} \\ 
-\frac{107}{943} & -\frac{1608}{943} & -\frac{620}{943} & \frac{1600}{943} & 
\frac{1111}{943} & \frac{469}{943} \\ 
-\frac{26}{205} & \frac{6}{41} & \frac{126}{205} & -\frac{99}{205} & -\frac{%
68}{205} & \frac{63}{205} \\ 
\frac{896}{4715} & \frac{629}{943} & \frac{89}{4715} & -\frac{1051}{4715} & -%
\frac{2482}{4715} & -\frac{878}{4715} \\ 
\frac{1147}{4715} & \frac{2025}{943} & \frac{4108}{4715} & -\frac{9202}{4715}
& -\frac{7124}{4715} & -\frac{3071}{4715} \\ 
-\frac{824}{4715} & -\frac{1412}{943} & -\frac{3576}{4715} & \frac{6734}{4715%
} & \frac{5903}{4715} & \frac{1902}{4715} \\ 
-\frac{507}{4715} & \frac{1388}{943} & \frac{4712}{4715} & -\frac{8388}{4715}
& -\frac{5426}{4715} & -\frac{924}{4715} \\ 
\frac{1188}{4715} & \frac{754}{943} & \frac{1402}{4715} & -\frac{2888}{4715}
& -\frac{2491}{4715} & -\frac{1964}{4715} \\ 
\frac{746}{4715} & -\frac{1276}{943} & -\frac{4041}{4715} & \frac{7934}{4715}
& \frac{4143}{4715} & \frac{132}{4715}%
\end{array}%
\right. $

$\ \ \ \ \ \ \ \ \ \ \ \ \ \ \ \ \ \ \ \ \ \ \ \ \ \ \ \ \ \ \ \ \ \ \ \
\left. 
\begin{array}{rrrrrr}
\frac{1562}{4715} & -\frac{1998}{4715} & -\frac{171}{943} & -\frac{689}{4715}
& \frac{1282}{4715} & \frac{3032}{4715} \\ 
-\frac{319}{943} & \frac{76}{943} & \frac{75}{943} & -\frac{105}{943} & -%
\frac{129}{943} & \frac{98}{943} \\ 
\frac{1056}{4715} & -\frac{89}{4715} & \frac{243}{943} & \frac{1128}{4715} & 
-\frac{1524}{4715} & \frac{456}{4715} \\ 
\frac{106}{205} & -\frac{139}{205} & -\frac{42}{41} & -\frac{157}{205} & 
\frac{271}{205} & \frac{111}{205} \\ 
-\frac{499}{943} & \frac{852}{943} & \frac{1486}{943} & \frac{1503}{943} & -%
\frac{2141}{943} & -\frac{837}{943} \\ 
\frac{78}{205} & -\frac{52}{205} & \frac{7}{41} & \frac{74}{205} & -\frac{52%
}{205} & \frac{43}{205} \\ 
-\frac{228}{4715} & -\frac{463}{4715} & -\frac{749}{943} & -\frac{3244}{4715}
& \frac{5687}{4715} & -\frac{1813}{4715} \\ 
\frac{3939}{4715} & -\frac{4471}{4715} & -\frac{2086}{943} & -\frac{8973}{%
4715} & \frac{13\,979}{4715} & \frac{3094}{4715} \\ 
-\frac{2653}{4715} & \frac{2657}{4715} & \frac{1455}{943} & \frac{5846}{4715}
& -\frac{7798}{4715} & -\frac{2753}{4715} \\ 
\frac{3981}{4715} & -\frac{3269}{4715} & -\frac{1278}{943} & -\frac{6142}{%
4715} & \frac{7596}{4715} & \frac{5041}{4715} \\ 
\frac{1971}{4715} & -\frac{1519}{4715} & -\frac{610}{943} & -\frac{4217}{4715%
} & \frac{5246}{4715} & \frac{101}{4715} \\ 
-\frac{3263}{4715} & \frac{5182}{4715} & \frac{1395}{943} & \frac{6266}{4715}
& -\frac{9168}{4715} & -\frac{4088}{4715}%
\end{array}%
\right] $

In addition the inverse of the anti-nonadiagonal matrix $Y$ is obtained as
follows:

\newpage

\item $Y^{-1}=RK^{-1}$

$=R\left[ 
\begin{array}{rrrrrr}
\frac{231}{4715} & \frac{199}{943} & \frac{3154}{4715} & -\frac{3181}{4715}
& -\frac{2187}{4715} & \frac{142}{4715} \\ 
\frac{79}{943} & -\frac{170}{943} & -\frac{591}{943} & \frac{643}{943} & 
\frac{440}{943} & -\frac{29}{943} \\ 
-\frac{1172}{4715} & \frac{15}{943} & \frac{3062}{4715} & -\frac{1088}{4715}
& -\frac{416}{4715} & \frac{96}{4715} \\ 
\frac{33}{205} & \frac{46}{41} & \frac{187}{205} & -\frac{308}{205} & -\frac{%
166}{205} & -\frac{9}{205} \\ 
-\frac{107}{943} & -\frac{1608}{943} & -\frac{620}{943} & \frac{1600}{943} & 
\frac{1111}{943} & \frac{469}{943} \\ 
-\frac{26}{205} & \frac{6}{41} & \frac{126}{205} & -\frac{99}{205} & -\frac{%
68}{205} & \frac{63}{205} \\ 
\frac{896}{4715} & \frac{629}{943} & \frac{89}{4715} & -\frac{1051}{4715} & -%
\frac{2482}{4715} & -\frac{878}{4715} \\ 
\frac{1147}{4715} & \frac{2025}{943} & \frac{4108}{4715} & -\frac{9202}{4715}
& -\frac{7124}{4715} & -\frac{3071}{4715} \\ 
-\frac{824}{4715} & -\frac{1412}{943} & -\frac{3576}{4715} & \frac{6734}{4715%
} & \frac{5903}{4715} & \frac{1902}{4715} \\ 
-\frac{507}{4715} & \frac{1388}{943} & \frac{4712}{4715} & -\frac{8388}{4715}
& -\frac{5426}{4715} & -\frac{924}{4715} \\ 
\frac{1188}{4715} & \frac{754}{943} & \frac{1402}{4715} & -\frac{2888}{4715}
& -\frac{2491}{4715} & -\frac{1964}{4715} \\ 
\frac{746}{4715} & -\frac{1276}{943} & -\frac{4041}{4715} & \frac{7934}{4715}
& \frac{4143}{4715} & \frac{132}{4715}%
\end{array}%
\right. $

$\ \ \ \ \ \ \ \ \ \ \ \ \ \ \ \ \ \ \ \ \ \ \ \ \ \ \ \ \ \ \ \ \ \ \ \ \ \
\ \ \ \ \ \left. 
\begin{array}{rrrrrr}
\frac{1562}{4715} & -\frac{1998}{4715} & -\frac{171}{943} & -\frac{689}{4715}
& \frac{1282}{4715} & \frac{3032}{4715} \\ 
-\frac{319}{943} & \frac{76}{943} & \frac{75}{943} & -\frac{105}{943} & -%
\frac{129}{943} & \frac{98}{943} \\ 
\frac{1056}{4715} & -\frac{89}{4715} & \frac{243}{943} & \frac{1128}{4715} & 
-\frac{1524}{4715} & \frac{456}{4715} \\ 
\frac{106}{205} & -\frac{139}{205} & -\frac{42}{41} & -\frac{157}{205} & 
\frac{271}{205} & \frac{111}{205} \\ 
-\frac{499}{943} & \frac{852}{943} & \frac{1486}{943} & \frac{1503}{943} & -%
\frac{2141}{943} & -\frac{837}{943} \\ 
\frac{78}{205} & -\frac{52}{205} & \frac{7}{41} & \frac{74}{205} & -\frac{52%
}{205} & \frac{43}{205} \\ 
-\frac{228}{4715} & -\frac{463}{4715} & -\frac{749}{943} & -\frac{3244}{4715}
& \frac{5687}{4715} & -\frac{1813}{4715} \\ 
\frac{3939}{4715} & -\frac{4471}{4715} & -\frac{2086}{943} & -\frac{8973}{%
4715} & \frac{13\,979}{4715} & \frac{3094}{4715} \\ 
-\frac{2653}{4715} & \frac{2657}{4715} & \frac{1455}{943} & \frac{5846}{4715}
& -\frac{7798}{4715} & -\frac{2753}{4715} \\ 
\frac{3981}{4715} & -\frac{3269}{4715} & -\frac{1278}{943} & -\frac{6142}{%
4715} & \frac{7596}{4715} & \frac{5041}{4715} \\ 
\frac{1971}{4715} & -\frac{1519}{4715} & -\frac{610}{943} & -\frac{4217}{4715%
} & \frac{5246}{4715} & \frac{101}{4715} \\ 
-\frac{3263}{4715} & \frac{5182}{4715} & \frac{1395}{943} & \frac{6266}{4715}
& -\frac{9168}{4715} & -\frac{4088}{4715}%
\end{array}%
\right] $

\newpage

\item $Y^{-1}=\left[ 
\begin{array}{cccccc}
\frac{746}{4715} & -\frac{1276}{943} & -\frac{4041}{4715} & \frac{7934}{4715}
& \frac{4143}{4715} & \frac{132}{4715} \\ 
\frac{1188}{4715} & \frac{754}{943} & \frac{1402}{4715} & -\frac{2888}{4715}
& -\frac{2491}{4715} & -\frac{1964}{4715} \\ 
-\frac{507}{4715} & \frac{1388}{943} & \frac{4712}{4715} & -\frac{8388}{4715}
& -\frac{5426}{4715} & -\frac{924}{4715} \\ 
-\frac{824}{4715} & -\frac{1412}{943} & -\frac{3576}{4715} & \frac{6734}{4715%
} & \frac{5903}{4715} & \frac{1902}{4715} \\ 
\frac{1147}{4715} & \frac{2025}{943} & \frac{4108}{4715} & -\frac{9202}{4715}
& -\frac{7124}{4715} & -\frac{3071}{4715} \\ 
\frac{896}{4715} & \frac{629}{943} & \frac{89}{4715} & -\frac{1051}{4715} & -%
\frac{2482}{4715} & -\frac{878}{4715} \\ 
-\frac{26}{205} & \frac{6}{41} & \frac{126}{205} & -\frac{99}{205} & -\frac{%
68}{205} & \frac{63}{205} \\ 
-\frac{107}{943} & -\frac{1608}{943} & -\frac{620}{943} & \frac{1600}{943} & 
\frac{1111}{943} & \frac{469}{943} \\ 
\frac{33}{205} & \frac{46}{41} & \frac{187}{205} & -\frac{308}{205} & -\frac{%
166}{205} & -\frac{9}{205} \\ 
-\frac{1172}{4715} & \frac{15}{943} & \frac{3062}{4715} & -\frac{1088}{4715}
& -\frac{416}{4715} & \frac{96}{4715} \\ 
\frac{79}{943} & -\frac{170}{943} & -\frac{591}{943} & \frac{643}{943} & 
\frac{440}{943} & -\frac{29}{943} \\ 
\frac{231}{4715} & \frac{199}{943} & \frac{3154}{4715} & -\frac{3181}{4715}
& -\frac{2187}{4715} & \frac{142}{4715}%
\end{array}%
\right. $

\ \ \ \ \ \ \ \ \ \ \ \ \ \ \ \ \ \ \ \ \ \ \ \ \ \ \ \ \ \ \ \ \ \ \ \ \ \
\ \ \ \ \ \ \ \ \ \ \ \ \ \ \ $\left. 
\begin{array}{cccccc}
-\frac{3263}{4715} & \frac{5182}{4715} & \frac{1395}{943} & \frac{6266}{4715}
& -\frac{9168}{4715} & -\frac{4088}{4715} \\ 
\frac{1971}{4715} & -\frac{1519}{4715} & -\frac{610}{943} & -\frac{4217}{4715%
} & \frac{5246}{4715} & \frac{101}{4715} \\ 
\frac{3981}{4715} & -\frac{3269}{4715} & -\frac{1278}{943} & -\frac{6142}{%
4715} & \frac{7596}{4715} & \frac{5041}{4715} \\ 
-\frac{2653}{4715} & \frac{2657}{4715} & \frac{1455}{943} & \frac{5846}{4715}
& -\frac{7798}{4715} & -\frac{2753}{4715} \\ 
\frac{3939}{4715} & -\frac{4471}{4715} & -\frac{2086}{943} & -\frac{8973}{%
4715} & \frac{13\,979}{4715} & \frac{3094}{4715} \\ 
-\frac{228}{4715} & -\frac{463}{4715} & -\frac{749}{943} & -\frac{3244}{4715}
& \frac{5687}{4715} & -\frac{1813}{4715} \\ 
\frac{78}{205} & -\frac{52}{205} & \frac{7}{41} & \frac{74}{205} & -\frac{52%
}{205} & \frac{43}{205} \\ 
-\frac{499}{943} & \frac{852}{943} & \frac{1486}{943} & \frac{1503}{943} & -%
\frac{2141}{943} & -\frac{837}{943} \\ 
\frac{106}{205} & -\frac{139}{205} & -\frac{42}{41} & -\frac{157}{205} & 
\frac{271}{205} & \frac{111}{205} \\ 
\frac{1056}{4715} & -\frac{89}{4715} & \frac{243}{943} & \frac{1128}{4715} & 
-\frac{1524}{4715} & \frac{456}{4715} \\ 
-\frac{319}{943} & \frac{76}{943} & \frac{75}{943} & -\frac{105}{943} & -%
\frac{129}{943} & \frac{98}{943} \\ 
\frac{1562}{4715} & -\frac{1998}{4715} & -\frac{171}{943} & -\frac{689}{4715}
& \frac{1282}{4715} & \frac{3032}{4715}%
\end{array}%
\right] $

where 
\[
R=\left[ 
\begin{array}{cccccccccccc}
0 & 0 & 0 & 0 & 0 & 0 & 0 & 0 & 0 & 0 & 0 & 1 \\ 
0 & 0 & 0 & 0 & 0 & 0 & 0 & 0 & 0 & 0 & 1 & 0 \\ 
0 & 0 & 0 & 0 & 0 & 0 & 0 & 0 & 0 & 1 & 0 & 0 \\ 
0 & 0 & 0 & 0 & 0 & 0 & 0 & 0 & 1 & 0 & 0 & 0 \\ 
0 & 0 & 0 & 0 & 0 & 0 & 0 & 1 & 0 & 0 & 0 & 0 \\ 
0 & 0 & 0 & 0 & 0 & 0 & 1 & 0 & 0 & 0 & 0 & 0 \\ 
0 & 0 & 0 & 0 & 0 & 1 & 0 & 0 & 0 & 0 & 0 & 0 \\ 
0 & 0 & 0 & 0 & 1 & 0 & 0 & 0 & 0 & 0 & 0 & 0 \\ 
0 & 0 & 0 & 1 & 0 & 0 & 0 & 0 & 0 & 0 & 0 & 0 \\ 
0 & 0 & 1 & 0 & 0 & 0 & 0 & 0 & 0 & 0 & 0 & 0 \\ 
0 & 1 & 0 & 0 & 0 & 0 & 0 & 0 & 0 & 0 & 0 & 0 \\ 
1 & 0 & 0 & 0 & 0 & 0 & 0 & 0 & 0 & 0 & 0 & 0%
\end{array}%
\right] .
\]
\end{itemize}


\begin{thebibliography}{9}
\bibitem{1} M. El-Mikkawy, E. D. Rahmo, Symbolic algorithm for
inverting\linebreak cyclic pentadiagonal matrices recursively - Derivation
and implementation, Computers \& Mathematics with Applications, 59 (2010)
Pages 1386-1396.

\bibitem{2} M.El-Mikkawy, E.D.Rahmo, A new recursive algorithm for inverting
general tridiagonal and anti-tridiagonal matrices, Appl. Math. Compt. 204
(2008) 368-372.

\bibitem{3} M.El-Mikkawy, E.D.Rahmo, A new recursive algorithm for
inverting\linebreak general periodic pentadiagonal and anti-pentadiagonal
matrices, Appl. Math. Compt. 207 (2009) 164-170.

\bibitem{4} M. El-Mikkawy, A fast algorithm for evaluating nth order
tri-diagonal\linebreak determinants, J. Compt. Appl. Math. 166 (2004)
581-584.

\bibitem{5} D. Aiat Hadj, M. Elouafi, A fast numerical algorithm for the
inverse of a tridiagonal and pentadiagonal matrix, Appl. Math. and Comp. 202
(2008) 441--445.

\bibitem{6} D. Aiat Hadj, M. Elouafi, On the characteristic polynomial,
eigenvectors and determinant of a pentadiagonal matrix, Appl. Math. and
Comp., 198 (2008), 634-642.
\end{thebibliography}
\end{document}